\documentclass[12pt,a4paper]{article}
\usepackage{amsmath,amsthm,amssymb,amscd}
\usepackage{graphicx,psfrag,comment}
\usepackage[dvips]{color}
\definecolor{Color1}{rgb}{1,0,0} 
\definecolor{Color2}{rgb}{0,0.5,0} 
\definecolor{Color3}{rgb}{0,0,1} 

\newcommand{\R}{\mathbb{R}}

\newcommand{\Z}{\mathbb{Z}}
\newcommand{\BF}{\mathbf{F}}
\newcommand{\D}{\mathcal{D}}
\newcommand{\E}{\mathcal{E}}
\newcommand{\X}{\mathcal{X}}
\newcommand{\Y}{\mathcal{Y}}
\newcommand{\tf}{\overline{f}}
\newcommand{\SP}{\operatorname{SP}}
\newcommand{\DSP}{\operatorname{DSP}}

 \DeclareMathOperator{\ind}{ind}
 \DeclareMathOperator{\tr}{tr}
 \DeclareMathOperator{\sgn}{sgn}
 \DeclareMathOperator{\intr}{int}
 \DeclareMathOperator{\fix}{Fix}
 \DeclareMathOperator{\cd}{cd}
 \DeclareMathOperator{\id}{id}

\newtheorem{thm}{Theorem}[section]
\newtheorem{lem}[thm]{Lemma}
\newtheorem{prop}[thm]{Proposition}

\theoremstyle{definition}
\newtheorem{defn}[thm]{Definition}
\newtheorem{exam}[thm]{Example}
\theoremstyle{remark}
\newtheorem{rem}[thm]{Remark}

\begin{document}

\title{A trace formula for the forcing relation of braids}
\author{
    Boju Jiang\footnote{E-mail: bjjiang@math.pku.edu.cn} , ~
    Hao Zheng\footnote{E-mail: zhenghao@mail.sysu.edu.cn}
    \medskip \\
    {\em \small Department of Mathematics, Peking University} \\
    {\em \small Beijing 100871, China}
}
\date{}
\maketitle

\begin{abstract}
The forcing relation of braids has been introduced for a
$2$-dimensional analogue of the Sharkovskii order on periods for
maps of the interval. In this paper, by making use of the Nielsen
fixed point theory and a representation of braid groups, we deduce
a trace formula for the computation of the forcing order.
\end{abstract}

\begin{small}

\smallskip

\end{small}

\section{Introduction}

The influential statement ``Period three implies chaos" of Li and
Yorke \cite{LY} turns out to be a consequence of a much earlier
theorem

\begin{thm}[Sharkovskii \cite{Sharkovskii}]
On the set of natural numbers, define a linear order
$$3 \succ 5 \succ 7 \succ \cdots \succ 2\cdot3 \succ 2\cdot5 \succ
  \cdots \succ 4\cdot3 \succ 4\cdot5 \succ \cdots
  \succ 8 \succ 4 \succ 2 \succ 1.
$$
For any continuous map $f : [0,1] \to [0,1]$ of the interval, if $f$
has a periodic orbit of period $n$, then $f$ must have a periodic
orbit of period $m$ for every $m \prec n$.
\end{thm}

However, in dimension $2$ the same statement cannot be true as shown
by the $(2\pi/3)$-rotation of the unit disk. It was not until the
work of Matsuoka \cite{Matsuoka} and Boyland \cite{Boyland1}, that
the role of braids was revealed in the problem of forcing relation
of periodic orbits for homeomorphisms of the plane.

\begin{figure}[h]
\centering
\includegraphics[scale=.6]{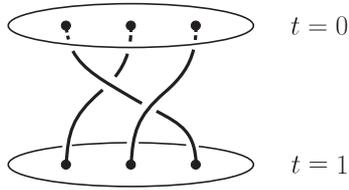}
\caption{ a geometric braid } \label{fig:fig11}
\end{figure}

Let $f : \R^2 \to \R^2$ be an orientation-preserving homeomorphism
and let $\{h_t : \R^2 \to \R^2\}_{0\le t\le1}$ be an isotopy with
$h_0 = \id$ and $h_1 = f$. An $f$-invariant set $P = \{
x_1,\dots,x_n \} \subset \R^2$ gives rise to a geometric braid (ref.
\cite{Birman})
$$\{ (h_t(x_i),t) \mid 0 \leq t \leq 1, \; 1 \leq i \leq n \}$$
in the cylinder $\R^2 \times [0,1]$. Indeed, the closed curve $\{
[h_t(x_1),\dots,h_t(x_n)] \mid 0 \leq t \leq 1\}$ in the
configuration space
$$\X_n = \{ (x_1,\dots,x_n) \mid x_i \in \R^2, \; x_i \neq x_j, \;
  \forall i \neq j \} / \Sigma_n,
$$
where $\Sigma_n$ denotes the symmetric group of $n$ symbols, gives
rise to a braid $\beta_P$ in the $n$-strand braid group $B_n =
\pi_1(\X_n)$. With another connecting isotopy $\{h_t\}$, the
resulting braid $\beta_P$ may differ by a power of the
``full-twist''. Matsuoka obtained lower bounds for the number of
$m$-periodic points of $f|_{\R^2 \setminus P}$, in terms of the
trace of the reduced Burau representation of the braid
$(\beta_P)^m$.

Later, Kolev \cite{Kolev} (see also \cite{Guaschi1}) found that a
$3$-periodic orbit $P$ guarantees the existence of $m$-periodic
orbits for every $m$, unless the braid $\beta_P$ is conjugate to a
power of the braid $\sigma_1\sigma_2$. Roughly speaking, this means
that the $(2\pi/3)$-rotation mentioned above is the only exceptional
case. Therefore, Li-Yorke's statement still holds in a subtle way in
$2$-dimensional dynamics. The analogue of the Sharkovskii order
naturally leads to the notion of forcing relation of (conjugacy
classes of) braids.

In the following, the notation $[\beta]$ stands for the conjugacy
class (in the group which is specified by the context) of a braid
$\beta$.

\begin{defn}\label{defn:forcing}
A braid $\beta$ {\em forces} a braid $\gamma$ if, for any
orientation-preserving homeomorphism $f : \R^2 \to \R^2$ and any
isotopy $\{h_t\} : \id \simeq f$, the existence of an $f$-invariant
set $P$ with $[\beta_P] = [\beta]$ guarantees the existence of an
$f$-invariant set $Q$ with $[\beta_{Q}] = [\gamma]$.
\end{defn}

\begin{rem}
There is a homomorphism from the braid group $B_n$ onto the mapping
class group of the pair $(\R^2,P)$ (acting on $(\R^2,P)$ \emph{from
the right}), its kernel being generated by the ``full-twist''. Via
this homomorphism, $[\beta_P]$ is sent to the conjugacy class of the
mapping class represented by $f$, which is independent of the choice
of the isotopy $\{h_t\}$. Following Boyland \cite{Boyland1} this
invariant is referred to as the {\em braid type} of $(f,P)$ in the
literature. It is clear that the forcing relation of braids defined
above naturally descends to that of braid types.
\end{rem}

The forcing relation is essentially a problem concerning plane
homeomorphisms. So the Bestvina-Handel theory of train-track maps
\cite{BH} comes in naturally. By analyzing the symbolic dynamics
of train-track maps, Handel \cite{Handel} was able to totally
solve the forcing relation among $3$-strand pseudo-Anosov braids,
and de Carvalho and Hall \cite{CH1, CH2} have managed to do the
same for horseshoe braids. This approach is, theoretical speaking,
powerful enough to be extended to mapping classes of all punctured
surfaces. On the other hand, there still is the challenging task
of recovering the braiding information encoded in the symbolic
dynamics.

\medskip

In this paper, we take another approach. Besides the Thurston
classification of surface homeomorphisms, we apply the Nielsen fixed
point theory. As a powerful tool for studying fixed points and
periodic orbits of self maps, Nielsen theory has been well developed
and successful in many mathematical problems. It turns out that
there are plenty of coincidences between the notions in the Nielsen
theory and the forcing theory, providing a more direct bridge
between the topological and the algebraic aspects of braids.

We start by slightly expanding the language of forcing.

\begin{defn}\label{defn:forced_extension}
A braid $\beta'$ is an \emph{extension} of $\beta$ if $\beta'$ is a
(disjoint but possibly intertwined) union of $\beta$ and another
braid $\gamma$. An extension $\beta'$ is \emph{forced by} $\beta$
if, for any orientation-preserving homeomorphism $f : \R^2 \to \R^2$
and any isotopy $\{h_t\} : \id \simeq f$, the existence of an
$f$-invariant set $P$ with $[\beta_P] = [\beta]$ guarantees the
existence of an additional $f$-invariant set $Q \subset \R^2
\setminus P$ with $[\beta_{P \cup Q}] = [\beta']$.
\end{defn}

The advantage of considering $[\beta_{P \cup Q}]$ is that it
contains the extra information of how the forced braid
$\gamma=\beta_Q$ winds around the original braid $\beta$.

\smallskip

Our main result is stated as follows.

\begin{thm}\label{thm:main}
Suppose a braid $\beta' \in B_{n+m}$ is an extension of $\beta \in
B_n$. Then $\beta'$ is forced by $\beta$ if and only if $\beta'$ is
neither collapsible nor peripheral relative to $\beta$, and the
conjugacy class $[\beta']$ has nonzero coefficient in
$\tr_{B_{n+m}}\zeta_{n,m}(\beta)$.
\end{thm}

In the theorem, $\zeta_{n,m}$ is a matrix representation of $B_n$
over a free $\Z B_{n+m}$-module, and the trace $\tr_{B_{n+m}}$ is
meant to take value in the free abelian group generated by the
\emph{conjugacy classes in $B_{n+m}$} (see Section
\ref{sec:formula}). In addition, $\beta'$ is said to be collapsible
or peripheral relative to $\beta$ if, roughly speaking, some strands
of $\beta'$ may be merged or moved to infinity while keeping $\beta$
untouched (see Definition \ref{defn:cp} and the figures therein).

Thus, to obtain the $(n+m)$-strand forced extensions of a braid
$\beta \in B_n$, it suffices to compute the trace
$\tr_{B_{n+m}}\zeta_{n,m}(\beta)$ and then drop off certain
irrelevant terms.

\medskip

The paper is organized as follows. In Section \ref{sec:Nielsen}, we
propose a Nielsen theory tailored to finite invariant sets of self
embeddings. In Section \ref{sec:forcing}, we apply this Nielsen
theory to the forcing problem of braids and reduce it to the
computation of a generalized Lefschetz number. In Section
\ref{sec:formula}, the representation $\zeta_{n,m}$ is defined and a
trace formula for the generalized Lefschetz number is derived.
Section \ref{sec:conjugacy} is devoted to the proof of Theorem
\ref{thm:main}. In the final section, we discuss the algorithmic
aspects of the trace formula and show some concrete examples.

The authors want to express their gratitude to the referee for
comments on the exposition, which motivated an extensive revision.

\section{Nielsen theory}\label{sec:Nielsen}

\subsection{Nielsen fixed point theory}\label{sec:Nielsen:orbit}

Throughout the paper all maps between topological spaces are assumed
to be continuous. The material in this subsection is standard, see
\cite{Jiang1, Jiang2}. We assume that $X$ is a compact, connected
polyhedron and $f : X \to X$ is a self map.

Consider the {\em mapping torus} $T_f = X \times \R_+ / (x,t+1) \sim
(f(x),t)$ of $f$. Denote by $\Gamma$ the fundamental group of $T_f$
and by $\Gamma_c$ the set of conjugacy classes of $\Gamma$. Then
$\Gamma_c$ is independent of the base point of $T_f$ and can be
regarded as the set of free homotopy classes of closed curves in
$T_f$.

Note that $x \in \fix f$ if and only if on the mapping torus $T_f$
its time-$1$ orbit curve $\{[x,t] \mid 0\le t\le 1 \}$ is closed.
Define $x,y \in \fix f$ to be in the same {\em fixed point class} if
and only if their time-$1$ orbit curves are freely homotopic in
$T_f$. Therefore every fixed point class $\BF$ gives rise to a
conjugacy class $\cd(\BF)$ in $\Gamma$, called the {\em coordinate}
of $\BF$. A fixed point class $\BF$ is called {\em essential}\/ if
its index $\ind(f,\BF)$ is nonzero.

The {\em generalized Lefschetz number} is defined as
$$L_\Gamma(f) = \sum_{\BF} \ind(f,\BF) \cdot \cd(\BF) \in \Z\Gamma_c$$
which takes value in the free abelian group $\Z\Gamma_c$ generated
by $\Gamma_c$.

The number of nonzero terms in $L_\Gamma(f)$ is called the {\em
Nielsen number} of $f$. It is the number of essential
fixed point classes, a lower bound for the number of fixed points of
$f$.

Generalized Lefschetz number is a homotopy invariant, i.e. if $f
\simeq g : X \to X$ then, identifying the fundamental groups of
$T_f$ and $T_g$ in the standard way, we have $L_\Gamma(f) =
L_\Gamma(g)$.

\subsection{Stratified maps}\label{sec:Nielsen:stratified}

The Nielsen theory for stratified maps is a version of relative
Nielsen theory. Readers are referred to \cite{JZZ} for a detailed
treatment of this subject.

\begin{defn}
Let $W$ be a compact, connected polyhedron and let $\emptyset=W^0
\subset W^1 \subset \dots \subset W^{m-1} \subset W^m = W$ be a
filtration of compact subpolyhedra. For $1\le k\le m$, the subspace
$W_k=W^k \setminus W^{k-1}$ is called the $k$-th stratum. A map $f :
W \to W$ is called a {\em stratified map} if $f(W_k) \subset W_k$
for all strata $W_k$. Two stratified maps $f,f' : W \to W$ are
called {\em stratified homotopic} if there is a homotopy of
stratified maps $\{h_t : W \to W\}_{0 \leq t \leq 1}$ such that $h_0
= f$, $h_1 = f'$.
\end{defn}

We will be concerned with fixed point classes of
$f_m=f|_{W_m}: W_m\to W_m$ in the top stratum.
A free homotopy class of closed curves in $T_{f_m}$,
represented by a closed curve $\gamma$, is said to be
{\em related} to a lower stratum $W_k$ if there is a homotopy of
closed curves $\{\gamma_s:[0,1]\to T_f\}_{0\le s\le1}$ such that
$\gamma_0=\gamma$, each $\gamma_s$
is in $T_{f_m}$ for $0\le s<1$, and
$\gamma_1$ is in $T_{f|_{W_k}}$.

\begin{defn}
A fixed point class of $f_m$ is called {\em degenerate} if its
coordinate is related to some lower stratum $W_k$. Otherwise, it is
called {\em non-degenerate}.
\end{defn}

Every non-degenerate fixed point class of $f_m$ is a compact subset
of $W_m$, hence its fixed point index is well defined.

Denote by $\Gamma$ the fundamental group of $T_{f_m}$ and by
$\Gamma_c$ the set of conjugacy classes in $\Gamma$. The {\em
generalized Lefschetz number} of the stratified map $f$ is defined
as
$$L_\Gamma(f)  = \sum_{\BF_m} \ind(f_m,\BF_m) \cdot \cd(\BF_m)
\in \Z\Gamma_c$$
where the sum is taken over all non-degenerate fixed point classes
$\BF_m$ of $f_m$.

The Nielsen fixed point theory has a natural version for stratified
maps. The main result is that $L_\Gamma(f)$ is not changed by a
stratified homotopy of the map $f$.

\smallskip

The following theorem will play an important role in the analysis of
the forcing relation of braids. Suppose $S$ is a compact, connected,
orientable surface and we consider homeomorphisms of $S$ as
stratified maps with respect to the filtration $\emptyset \subset
\partial S \subset S$.

\begin{thm}[Jiang-Guo \cite{JG}, Boyland \cite{Boyland2}]\label{thm:min2}
Every orientation-preserving homeomorphism $f : S \to S$ is isotopic
(through homeomorphisms) to a homeomorphism $\phi$ such that, for
any $n\ge1$, any fixed point class of $\phi^n$ is essential and
contained in a single $\phi$-orbit and, moreover, no fixed point
class of $\phi^n|_{\intr S}$ is related to $\partial{S}$.
\end{thm}

In the theorem, the $\phi|_{\intr S}$-orbits persist under isotopy,
i.e. none of them can be merged or be eliminated by isotoping the
homeomorphism $\phi|_{\intr S}$. In particular, $\phi|_{\intr S}$
has the minimal number of periodic orbits of period $n$ in its
isotopy class for every $n\ge1$. In the rest of this paper, we will
refer to the homeomorphism $\phi|_{\intr S}$ as a {\em minimal
representative} in the isotopy class of $f|_{\intr S}$.

\subsection{A Nielsen theory for finite invariant sets}
\label{sec:Nielsen:set}

In this subsection, assume $X$ is a compact, connected, smooth
manifold of dimension $d$ and $f : X \to X$ is a self
\emph{embedding}.

Let $m$ be a fixed natural number. Consider the symmetric product
space
$$\SP^m X = \{ (x_1,\dots,x_m) \mid x_i \in X \} / \Sigma_m. $$
Its points will be written as $[x_1,\dots,x_m]$, with repetitions
allowed. For an integer $k$, $0\le k \le m$, define the subspace
$$\SP^{m,k} X = \{ [x_1,\dots,x_m] \in \SP^m X \mid |\{x_1,\dots,x_m\}| \leq k \}.$$
Then we have a filtration
$$\emptyset=\SP^{m,0} X \subset \SP^{m,1} X \subset \dots \subset
\SP^{m,m-1} X \subset \SP^{m,m} X = \SP^m X.$$
The top stratum is sometimes called the deleted $m$-th symmetric
product space and denoted
$$\DSP^m X= \SP^m X \setminus \SP^{m,m-1} X.$$
For $1\le k\le m$, the $k$-th stratum is $W_k=\SP^{m,k} X \setminus
\SP^{m,k-1} X$.

\begin{rem}
Each stratum $W_k$ is a manifold of dimension $k\cdot d$. When
$d=2$, which is our main concern later, the $m$-th symmetric product
$\SP^m X$ itself is a manifold of dimension $2m$.
\end{rem}

The map $f:X\to X$ induces a map $\SP^m f : \SP^m X \to \SP^m X$
given by $\SP^m f([x_1,\dots,x_m])=[f(x_1),\dots,f(x_m)]$. Since $f$
is an embedding, $\SP^m f$ is now a stratified map with respect to
the above filtration. Hence the theory in the previous subsection is
applicable.

Observe that a fixed point $[x_1,\dots,x_m]$ of $\DSP^m f= \SP^m
f|_{\DSP^m X}$ corresponds to an $f$-invariant set consisting of
precisely $m$ distinct points. Thus, the number of non-degenerate,
essential fixed point classes of $\DSP^m f$ is a lower bound for the
number of such $f$-invariant sets, for all embeddings isotopic to
$f$.

Below is a useful criterion for the degeneracy of a fixed point
class of $\DSP^mf$.

\begin{prop}\label{prop:degenerate}
Suppose $X$ is a compact, connected, smooth manifold of dimension
$d$ and suppose $f:X\to X$ is a self embedding. Let $Q =\{
x_1,\dots,x_m \} \subset X$ be an $f$-invariant set. Let $\D$ denote
the disjoint union of $k$ copies of the $d$-dimensional disk, $1\le
k<m$. The coordinate of the fixed point $[x_1,\dots,x_m]$ of $\DSP^m
f$ is related to the $k$-th stratum $W_k$ if and only if there
exists an isotopy of embeddings $\{i_t : \D \to X\}_{0\le t\le1}$
such that $i_0 = f \circ i_1$, $Q \subset i_t(\D)$ and each
component of $i_t(\D)$ contains at least one point of $Q$, for all
$0 \leq t \leq 1$.
\end{prop}
\begin{proof}
Sufficiency is clear: we can use the tubes $\{ (x,t) \mid x\in
i_t(\D), \; t\in[0,1] \}$ to construct a homotopy which relates the
time-$1$ orbit curve of $[x_1,\dots,x_m]$ to a closed curve in the
mapping torus of $\SP^m f|_{W_k}$.

Necessity. Let $\{\gamma_s : [0,1] \to T_{\SP^m f}\}_{s\in[0,1]}$
be a homotopy of closed curves, relating the time-$1$ orbit curve
$\gamma$ of $[x_1,\dots,x_m]$ to a closed curve lying in
the mapping torus of $\SP^m f|_{W_k}$. By sliding along
$T_{\SP^m f}$, we may assume
that $\gamma_s(t)=[[z_1(s,t),\dots,z_m(s,t)],t]$,
where $z_j(s,t)$ are continuous functions of $s,t$.
Put $Q_{s,t} =\{ z_j(s,t) \mid 1 \leq j \leq m \}\subset X$. Then
$Q_{0,t} = Q$ and $Q_{s,0} = f(Q_{s,1})$.

Since there are precisely $k$ distinct elements in $Q_{1,t}$, if
$z_{j_1}(1,t)=z_{j_2}(1,t)$ holds for some $t$ then it holds for all
$t$. Therefore, we may choose an isotopy of embeddings $\{i'_t: \D \to
X\}_{t\in[0,1]}$ such that $i'_0 = f \circ i'_1$ and each component of
$i'_t(\D)$ contains precisely one point of $Q_{1,t}$. Extending
$\{i'_t\}$ to a two-parameter isotopy $\{i''_{s,t} : \D \to X\}_{s,t\in[0,1]}$ such that
$i''_{1,t} = i'_t$, $i''_{s,0} = f \circ i''_{s,1}$ and
$Q_{s,t} \subset i''_{s,t}(\D)$, we get the desired isotopy $\{i_t =
i''_{0,t}: \D \to X\}_{t\in[0,1]}$.
\end{proof}

\begin{defn}
In Proposition \ref{prop:degenerate}, the components of $i_0(\D)$
containing more than one point of $Q$ are called {\em merging disks}
of $Q$.
\end{defn}

\begin{rem}\label{rem:degenerate}
The existence of merging disks of $Q$ means the $f$-invariant set
$Q$ can be merged into a smaller one by isotoping $f$ in a
neighborhood of these disks.
\end{rem}

\subsection{Index formulae}\label{sec:Nielsen:index}

The next two lemmas may be found in \cite{Jiang1}.

\begin{lem}\label{lem:ind1}
Suppose $x$ is a generic fixed point of $f : X \to X$, i.e. $f$ is
differentiable at $x$ with Jacobian $A$ such that $\det(I-A) \neq
0$. Then $x$ is an isolated fixed point and $\ind(f,x) = \sgn
\det(I-A)$.
\end{lem}

\begin{lem}\label{lem:ind2}
Suppose $x$ and $y$ are isolated fixed points of $f : X \to X$ and
$g : Y \to Y$, respectively. Then $\ind(f \times g, x \times y) =
\ind(f,x) \cdot \ind(g,y)$.
\end{lem}

\begin{lem}\label{lem:ind3}
Let $f : (\R^n)^k \to (\R^n)^k$ be a map defined by
$$(x_1,\dots,x_k) \mapsto (f_2(x_2),\dots,f_{k}(x_{k}),f_1(x_1))$$
where $f_1,\dots,f_k : \R^n \to \R^n$ are a family of maps, each
admitting the origin $0\in\R^n$ as a fixed point. If the origin
$0\in(\R^n)^k $ is an isolated fixed point of $f$, then $\ind(f,0) =
\ind(f_1 \circ\cdots\circ f_k, 0)$.
\end{lem}
\begin{proof}
It suffices to consider the case of generic fixed points. Denote
by $A_1,\cdots,A_k$ the Jacobians of $f_1,\cdots,f_k$ at $0$. Then
the Jacobian of $f$ at $0$ is
$$A = \left( \begin{array}{ccccc}
    0   & A_2                           \\
        & 0     & A_3                   \\
        &       & 0                     \\
        &       &       & \ddots & A_k  \\
    A_1 &       &       &       &  0
  \end{array} \right).
$$
Therefore,
$$\ind(f,0) = \sgn \det(I-A) = \sgn \det(I-A_1 \cdots A_k)
= \ind(f_1 \circ\cdots\circ f_k, 0).
\eqno\qed
$$
\renewcommand\qed{}
\end{proof}

The following lemma is prepared for the proof of Proposition
\ref{prop:trace}.

\begin{lem}\label{lem:ind5}
Let $\lambda > 1$ be a real number and let $B=(B_{ij})$ be an $m
\times m$ matrix with
$$B_{ij} = \left\{ \begin{array}{ll}
    (-1)^{n_j}, & \eta(j) = i, \\
    0, & \eta(j) \neq i,
  \end{array} \right.
$$
where $n_1,\dots,n_m \in \Z$ are a set of integers and $\eta \in
\Sigma_m$ is a permutation. For a generic fixed point $x$ of a map
$f : \R^{2m} \to \R^{2m}$ with Jacobian $A =
\begin{pmatrix} \lambda B \\ & \lambda^{-1}B \end{pmatrix}$, we have
$$\ind(f,x) = (-1)^m (-1)^{n_1+\cdots+n_m} \sgn\eta.$$
\end{lem}
\begin{proof}
We may assume $\eta$ is a cycle. Then
$$\det(I-A) = (1-(-1)^{n_1+\cdots+n_m}\lambda^m)
  (1-(-1)^{n_1+\cdots+n_m}\lambda^{-m}).
$$
Hence
$$\ind(f,x) = \sgn \det(I-A)
  = -(-1)^{n_1+\cdots+n_m}
  = (-1)^m (-1)^{n_1+\cdots+n_m} \sgn\eta.
\eqno\qed
$$
\renewcommand\qed{}
\end{proof}

\section{Forced extensions of a braid}\label{sec:forcing}

In this section, we apply Nielsen theory to the problem of forcing
relation of braids. It turns out that the mapping tori involved here
are naturally embedded into the defining configuration spaces of
braid groups, and the coordinates of the fixed points can be readily
interpreted as forced extensions.

\subsection{Coordinates recognized as braid extensions}

Let $\sigma_1,\dots,\sigma_{n-1}$ denote the standard generators of
the Artin's $n$-strand braid group $B_n$ (ref. \cite{Birman}) and
let $A_{i,j}$, $1\le i<j\le n$, denote the standard pure braid
$$A_{i,j} = \sigma_{j-1} \cdots \sigma_{i+1} \sigma_i^2
  \sigma_{i+1}^{-1} \cdots \sigma_{j-1}^{-1}.
$$
\begin{figure}[h]
\centering
\includegraphics[scale=.6]{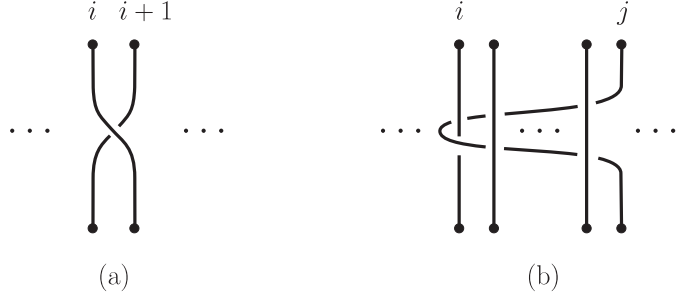}
\caption{ the standard braids (a) $\sigma_i$ and (b) $A_{i,j}$ }
\label{fig:fig13}
\end{figure}

Let $P \subset \R^2$ be a prescribed set of $n$ punctures. Define
the configuration spaces
\begin{eqnarray*}
  && \X_n = \{ (x_1,\dots,x_n) \mid x_i \in \R^2, \; x_i \neq x_j, \;
    \forall i \neq j \} / \Sigma_n, \\
  && \X_{n,m} = \{ (x_1,\dots,x_{n+m}) \mid x_i \in \R^2, \; x_i \neq x_j, \;
    \forall i \neq j \} / \Sigma_n \times \Sigma_m, \\
  && \Y_{n,m} = \{ (y_1,\dots,y_m) \mid y_i \in \R^2 \setminus P, \;
    y_i \neq y_j, \; \forall i \neq j  \} / \Sigma_m.
\end{eqnarray*}
Clearly $\Y_{n,m}$ embeds into $\X_{n,m}$ via
$[y_1,\dots,y_m]\mapsto[P,y_1,\dots,y_m]$, and $\X_{n,m}$ projects
onto $\X_n$ via
$[x_1,\dots,x_n,x_{n+1},\dots,x_{n+m}]\mapsto[x_1,\dots,x_n]$. Then
$\Y_{n,m}$ is precisely the fiber of the bundle $\pi : \X_{n,m} \to
\X_n$ and we have
\begin{eqnarray*}
  && \pi_1(\X_n) = B_n, \\
  && \pi_1(\X_{n,m}) = \langle \sigma_1,\dots,\sigma_{n-1},\sigma_n^2,
    \sigma_{n+1},\dots,\sigma_{n+m-1} \rangle \subset B_{n+m}, \\
  && \pi_1(\Y_{n,m}) = \langle A_{1,n+1},\dots,A_{n,n+1},\sigma_{n+1},\dots,
    \sigma_{n+m-1} \rangle \subset \pi_1(\X_{n,m}).
\end{eqnarray*}

Given a nontrivial $n$-strand braid $\beta$, ``sliding the plane
down the braid $\beta$'' determines (up to isotopy) a homeomorphism
$f_\beta : \R^2 \setminus P \to \R^2 \setminus P$, as well as a
connecting isotopy $\{h_t: \R^2 \to \R^2\}_{0\le t\le1}:\id \simeq
f_\beta$ such that the curves $\{h_t(P)\}_{0\le t\le1}$ represent
the braid $\beta$.

Now we figure out our key observations.
\begin{prop}\label{prop:key}
(1) The mapping torus of the induced map $\DSP^mf_\beta:
\Y_{n,m}\to\Y_{n,m}$ can be identified with the space obtained from
$$
\{ ([h_t(P),y_1,\dots,y_m],t) \mid y_i \in \R^2 \setminus h_t(P), \;
0\le t\le1 \} \subset \X_{n,m} \times [0,1]
$$
by identifying the top $\Y_{n,m} \times 0$ with the bottom $\Y_{n,m}
\times 1$.

(2) Under the above identification, the fundamental group
$\Gamma_{\beta,m}$ of\/ $T_{\DSP^mf_\beta}$ is the subgroup in
$B_{n+m}$ generated by $\beta$ and $\pi_1(\Y_{n,m})$, where $\beta$
is regarded as an $(n+m)$-strand braid with $m$ trivial strands
added.

(3) Moreover, when a fixed point of $\DSP^mf_\beta$ corresponds to
an $f_\beta$-invariant set $Q\subset\R^2\setminus P$, the coordinate
of the former is precisely $[\beta_{P\cup Q}]$.
\end{prop}

\begin{proof}
(1) The identification is realized by the embedding
\begin{align*}
  \Y_{n,m}\times[0,1] &\to \X_{n,m}\times[0,1], \\
  ([P,y_1,\dots,y_m],t) &\mapsto ([h_t(P),h_t(y_1),\dots,h_t(y_m)],t).
\end{align*}
(2) and (3) follow from (1).
\end{proof}

\subsection{Compactification}

By Proposition \ref{prop:key}, the coordinates of the fixed points
of $\DSP^mf_\beta$ are naturally interpreted as $(n+m)$-strand
extensions of $\beta$. However, $\DSP^mf_\beta$ is a self map of a
noncompact space, hence falls out of the framework of the usual
Nielsen fixed point theory.

So, we apply instead the theory in Section \ref{sec:Nielsen:set}.
For this, we compactify $\R^2 \setminus P$ to a $2$-disk with $n$
holes and denote it by $Y_n$, and assume further that $f_\beta =
\tf_\beta|_{\intr Y_n}$ for some homeomorphism $\tf_\beta : Y_n \to
Y_n$.

Consider the symmetric product space $\SP^m Y_n$ and the induced
stratified map $\SP^m\tf_\beta: \SP^m Y_n\to\SP^m Y_n$ with respect
to the filtration $\SP^{m,0} Y_n \subset \SP^{m,1} Y_n \subset \dots
\subset \SP^{m,m} Y_n$. Note that the interior of the manifold
$\DSP^m Y_n$ is precisely $\Y_{n,m}$. In particular, $\pi_1(\DSP^m
Y_n) = \pi_1(\Y_{n,m})$ and $\pi_1(T_{\DSP^m\tf_\beta}) =
\pi_1(T_{\DSP^mf_\beta}) = \Gamma_{\beta,m}$.

The stratified map $\SP^m\tf_\beta$ is actually the desired
compactification of $\DSP^mf_\beta$. But some more fixed points may
arise on the boundary of the manifold $\DSP^m Y_n$. Hence the
generalized Lefschetz number $L_{\Gamma_{\beta,m}}(\SP^m\tf_\beta)$
may contain unwanted terms (called ``peripheral'' terms) which
should be identified and ruled out.

In addition, the coordinates of degenerate fixed point classes of
$\DSP^m\tf_\beta$ also need to be identified in the computation of
$L_{\Gamma_{\beta,m}}(\SP^m\tf_\beta)$. These considerations lead to
the notions of the next subsection.

\subsection{Collapsible and peripheral extensions}

Recall the Thurston classification theorem for homeomorphisms of
compact surfaces. The theorem has a natural version for punctured
surfaces, even for punctured planes by regarding the plane as a
once-punctured $2$-sphere.

\begin{thm}[Thurston \cite{FLP, Thurston}]
Every homeomorphism $f : S \to S$ of a compact surface $S$ is
isotopic to a homeomorphism $\phi$ (Thurston representative) such
that either

(1) $\phi$ is a periodic map, i.e. $\phi^k = \id$ for some positive
integer $k$; or

(2) $\phi$ is a pseudo-Anosov map, i.e. there is a number
$\lambda>1$ and a pair of transverse measured foliations
$(F^s,\mu^s)$ and $(F^u,\mu^u)$ such that $\phi(F^s,\mu^s) =
(F^s,\lambda^{-1}\mu^s)$ and $\phi(F^u,\mu^u) = (F^u,\lambda\mu^u)$;
or

(3) $\phi$ is a reducible map, i.e. there is a system of disjoint
simple closed curves $\gamma = \{\gamma_1,\dots,\gamma_k\}$ in
$\intr S$ (reducing curves) such that $\gamma$ is invariant by
$\phi$ (but $\gamma_i$'s may be permuted) and $\gamma$ has a
$\phi$-invariant tubular neighborhood $U$ such that each component
of $S \setminus U$ has negative Euler characteristic and on each
$\phi$-component of $S \setminus U$, $\phi$ satisfies (1) or (2).
\end{thm}

Every braid determines a unique isotopy class of homeomorphisms of a
punctured plane. In this way, the braids naturally fall into three
types: periodic, pseudo-Anosov and reducible.

\begin{defn}\label{defn:cp}
Suppose $\beta' \in B_{n+m}$ is an extension of $\beta \in B_n$. Let
$\phi$ be a Thurston representative determined by $\beta'$. We say
$\beta'$ is {\em collapsible} (resp.\ {\em peripheral}) relative to
$\beta$ if there exists a system of reducing curves of $\phi$ such
that one of them encloses none of (resp.\ precisely one of or all
of) the punctures corresponding to $\beta$.
\end{defn}

\begin{figure}[h]
\centering
\includegraphics[scale=.6]{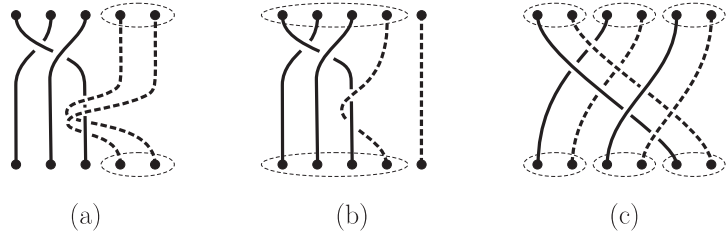}
\caption{ (a) collapsible and (b,c) peripheral braids relative to
the solid braid } \label{fig:fig12}
\end{figure}

\begin{defn}
If an extension $\beta' \in \beta\cdot\pi_1(\Y_{n,m})$ of a braid
$\beta\in B_n$ is collapsible (resp.\ peripheral) relative to
$\beta$, then we say the conjugacy class $[\beta']$ in
$\Gamma_{\beta,m}$ is {\em collapsible} (resp.\ {\em peripheral}).
\end{defn}

\subsection{Forced extensions}

We are ready to state the main result of this section.

\begin{prop}\label{prop:Lefschetz}
Suppose $\beta \in B_n$ is a nontrivial braid. The $(n+m)$-strand
forced extensions of $\beta$ are exactly the non-peripheral terms in
$L_{\Gamma_{\beta,m}}(\SP^m\tf_\beta)$.
\end{prop}

\begin{proof}
Thanks to the homotopy invariance of generalized Lefschetz number,
let us assume $f_\beta=\tf_\beta|_{\intr Y_n}: \R^2\setminus
P\to\R^2\setminus P$ is a minimal representative (in the sense of
Theorem \ref{thm:min2}) in its isotopy class.

On the one hand, a fixed point $[x_1,\dots,x_m]$ of
$\DSP^m\tf_\beta$ falls out of $\Y_{n,m}$ if and only if some $x_i$
falls into $\partial{Y}_n$, and this is equivalent by the minimality
of $f_\beta$ to that the coordinate of the fixed point is
peripheral.

On the other hand, the fixed point class represented by a fixed
point $[x_1,\dots,x_m]$ of $\DSP^m\tf_\beta$ lying in $\Y_{n,m}$ is
non-degenerate; otherwise, by Proposition \ref{prop:degenerate} the
$f_\beta$-invariant set $\{x_1,\dots,x_m\}$ can be merged into a
smaller one (cf. Remark \ref{rem:degenerate}), in contradiction to
the minimality of $f_\beta$.

It follows that the non-peripheral terms in
$L_{\Gamma_{\beta,m}}(\SP^m\tf_\beta)$ are precisely the coordinates
of the fixed points of $\DSP^mf_\beta =
\DSP^m\tf_\beta|_{\Y_{n,m}}$, which by the minimality of $f_\beta$
again are exactly the $(n+m)$-strand forced extensions of $\beta$.
\end{proof}

From the theoretical point of view, constructing a minimal
representative in the isotopy class of $f_\beta$ is not an easy
task. However, by the homotopy invariance of generalized Lefschetz
number, we can compute $L_{\Gamma_{\beta,m}}(\SP^m\tf_\beta)$ from
any map stratified homotopic to $\SP^m\tf_\beta$. This is exactly
what we will do in the next section. The following lemma is prepared
for this purpose.

\begin{lem}\label{lem:collapse}
Let $g : Y_n \to Y_n$ be an embedding isotopic to $\tf_\beta$ and
let $Q = \{ x_1,\dots,x_m \}$ be a $g$-invariant set. Then the fixed
point class of $\DSP^mg$ represented by $[x_1,\dots,x_m]$ is
degenerate if and only if its coordinate is collapsible.
\end{lem}

\begin{proof}
There is an obvious correspondence between the merging disks of a
$g$-invariant set and the reducing curves defining the notion of
collapsibility.
\end{proof}

\section{A trace formula for $L_{\Gamma_{\beta,m}}(\SP^m\tf_\beta)$}\label{sec:formula}

In this section, we define a representation $\zeta_{n,m}$ of $B_n$
over the free (left) $\Z B_{n+m}$-module generated by
$$\E_{n,m} = \{ \mu=(\mu_1,\dots,\mu_{n-1}) \mid \mu_i \in \Z_{\geq 0},
  \; \mu_1 + \cdots + \mu_{n-1} = m \},
$$
and derive a trace formula for
$L_{\Gamma_{\beta,m}}(\SP^m\tf_\beta)$. The size of the basis
$\E_{n,m}$ is $\binom{m+n-2}{m}$.

\subsection{Braid actions on $Y_n$}

\begin{figure}[h]
\centering
\includegraphics[scale=.6]{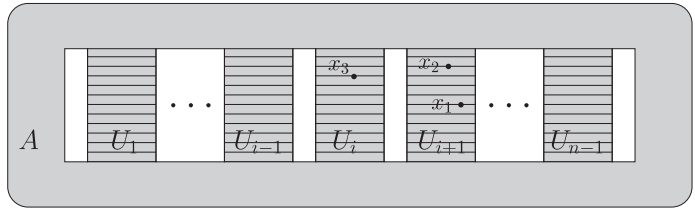}
\caption{ decomposition of $Y_n$ } \label{fig:fig21}
\end{figure}

We decompose the surface $Y_n$ into an annulus and $n-1$ rectangles,
as shown in Figure \ref{fig:fig21}. Let $U = U_1 \cup \cdots \cup
U_{n-1}$ be the union of the $n-1$ foliated open rectangles. Define
an ordering on $U$ such that $x_1 \prec x_2$ if either $x_1$ lies in
a rectangle to the right of $x_2$ or $x_1$ lies in a strictly lower
leaf of the same rectangle as $x_2$. For example, the order of the
three points in Figure \ref{fig:fig21} is $x_1 \prec x_2 \prec x_3$.

Set
$$V = \{ [x_1,\dots,x_m] \in \Y_{n,m} \mid x_i \in U, \; x_1 \prec \cdots \prec x_m \}.$$
Then $V = \bigcup_{\mu \in \E_{n,m}} V_\mu$ where
$$V_\mu = \{ [x_1,\dots,x_m] \in V \mid |\{x_1,\dots,x_m\} \cap U_i| = \mu_i,
  \; 1 \leq i \leq n-1 \}.
$$
Each $V_\mu$ is connected, thus the elements of $\E_{n,m}$ are in
one-one correspondence to the components of $V$.

\begin{figure}[h]
\centering
\includegraphics[scale=.6]{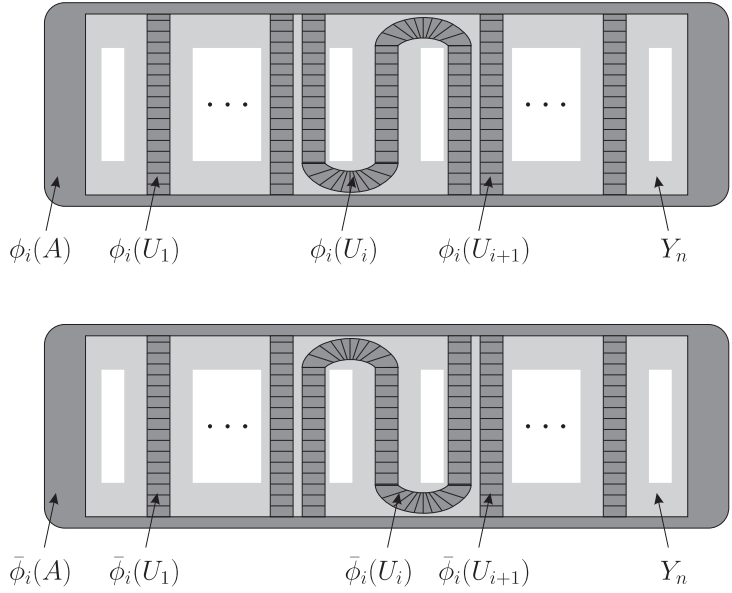}
\caption{ the embeddings $\phi_i : Y_n \to Y_n$ and $\bar\phi_i :
Y_n \to Y_n$ } \label{fig:fig22}
\end{figure}

Illustrated in Figure \ref{fig:fig22} are two embeddings $\phi_i :
Y_n \to Y_n$ and $\bar\phi_i : Y_n \to Y_n$, which can be understood
as the action of the elementary braids $\sigma_i$ and$\sigma_i^{-1}$
on $Y_n$, respectively. Both push the annulus outward, irrationally
rotate the outmost boundary, keep the foliations of $\phi_i^{-1}(U)$
or $\bar\phi_i^{-1}(U)$, uniformly contract along the leaves of the
foliations and uniformly expand along the transversal direction.
Slightly abusing notations, we also use $\phi_i$ and $\bar\phi_i$ to
denote the induced stratified maps of $\SP^mY_n$.

For every $\phi \in \{ \phi_1,\dots,\phi_{n-1},
\bar\phi_1,\dots,\bar\phi_{n-1} \}$, we have
$$V_\mu \cap \phi^{-1}(V_\nu) = \bigcup_{\eta \in \Sigma_m} W^{(\phi)}_{\mu\nu\eta}$$
where
$$W^{(\phi)}_{\mu\nu\eta} = \bigg\{ [x_1,\dots,x_m] \in V_\mu \cap \phi^{-1}(V_\nu)
  \bigg | \begin{array}{l}
    x_{\eta(1)} \prec \cdots \prec x_{\eta(m)}, \\
    \phi(x_1) \prec \cdots \prec \phi(x_m)
  \end{array} \bigg \}.
$$
Each $W^{(\phi)}_{\mu\nu\eta}$ is connected, thus the elements of
the set $\{ \eta \in \Sigma_m \mid W^{(\phi)}_{\mu\nu\eta} \neq
\emptyset \}$ are in one-one correspondence to the components of
$V_\mu \cap \phi^{-1}(V_\nu)$.

\begin{figure}[h]
\centering
\includegraphics[scale=.6]{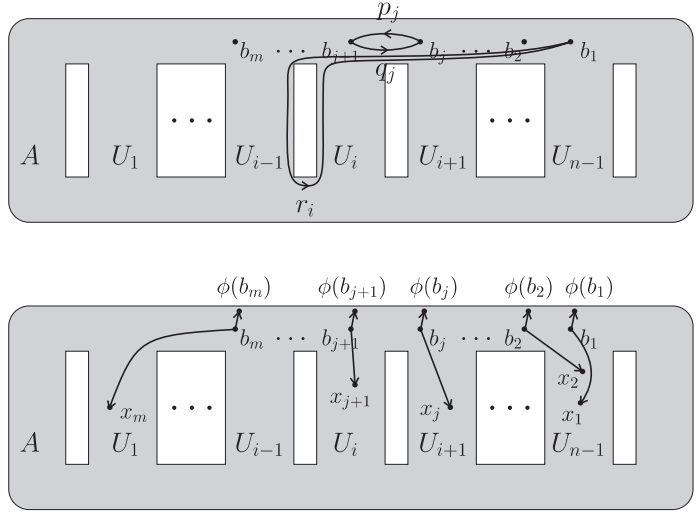}
\caption{ paths in $Y_n$ } \label{fig:fig23}
\end{figure}

Further, as shown in Figure \ref{fig:fig23}, choose a base point $b
= [b_1,\dots,b_m] \in \Y_{n,m}$. Then the generators $\sigma_{n+j}$
and $A_{i,n+1}$ of $\pi_1(\Y_{n,m})$ are represented by the loops
$[b_1,\dots,b_{j-1},p_j,q_j,b_{j+2},\dots,b_m]$ and
$[r_i,b_2,\dots,b_m]$, respectively.

For every $x=[x_1,\dots,x_m] \in V$ with $x_1 \prec \cdots \prec
x_m$, the disjoint ``descending" paths connecting $b_k$ to $x_k$ in
$Y_n$ give rise to a path $\gamma_x$ in $\Y_{n,m}$. Similarly, the
disjoint ``ascending" paths connecting $b_k$ to $\phi(b_k)$ give
rise to a path $\gamma_{\phi(b)}$ in $\Y_{n,m}$. For every nonempty
$W^{(\phi)}_{\mu\nu\eta}$, we choose a point $x \in
W^{(\phi)}_{\mu\nu\eta}$ and let $\alpha^{(\phi)}_{\mu\nu\eta}$
denote the element of $\pi_1(\Y_{n,m})$ represented by the loop
$\gamma_{\phi(b)} \cdot \phi(\gamma_x) \cdot \gamma^{-1}_{\phi(x)}$.
Note that $\alpha^{(\phi)}_{\mu\nu\eta}$ is independent of the
choices of $x$, $\gamma_x$, $\gamma_{\phi(b)}$ and
$\gamma_{\phi(x)}$.

\subsection{The representation $\zeta_{n,m}$}

\begin{prop}
The equations
\begin{eqnarray*}
  && \mu \cdot \zeta_{n,m}(\sigma_i)
  = \sum_{\nu \in \E_{n,m}} c^{(i)}_{\mu\nu} \cdot \nu, \\
  && \mu \cdot \zeta_{n,m}(\sigma_i^{-1})
  = \sum_{\nu \in \E_{n,m}} d^{(i)}_{\mu\nu} \cdot \nu,
\end{eqnarray*}
where
\begin{eqnarray*}
  c^{(i)}_{\mu\nu} &=& (-1)^{\nu_i} \cdot \sigma_i \cdot
  \sum_{\eta : W^{(\phi_i)}_{\mu\nu\eta} \neq \emptyset}
  \sgn\eta \cdot \alpha^{(\phi_i)}_{\mu\nu\eta}, \\
  d^{(i)}_{\mu\nu} &=& (-1)^{\nu_i} \cdot \sigma_i^{-1} \cdot
  \sum_{\eta : W^{(\bar\phi_i)}_{\mu\nu\eta} \neq \emptyset}
  \sgn\eta \cdot \alpha^{(\bar\phi_i)}_{\mu\nu\eta},
\end{eqnarray*}
give rise to a group representation $\zeta_{n,m}$ of $B_n$ over
the free $\Z B_{n+m}$-module generated by $\E_{n,m}$.
\end{prop}
\begin{proof}
Direct proof of this proposition is rather complicated. We refer the
reader to \cite{Zheng1} for another representation $\xi_{n,m}$ of
$B_n$ with the same basis $\E_{n,m}$ and with the same coefficient
ring $\Z B_{n+m}$. Let $a^{(i)}_{\mu\nu}$, $b^{(i)}_{\mu\nu}$ denote
the matrix elements of $\xi_{n,m}(\sigma_i)$,
$\xi_{n,m}(\sigma_i^{-1})$, respectively. Then a direct computation
shows that $c^{(i)}_{\mu\nu} = (b^{(i)}_{\nu\mu})^*$ and
$d^{(i)}_{\mu\nu} = (a^{(i)}_{\nu\mu})^*$, where $* : \Z B_{n+m} \to
\Z B_{n+m}$ is the involution determined by $a^* = a^{-1}$ for $a
\in B_{n+m}$. Therefore, $\zeta_{n,m}$ is precisely the dual
representation of $\xi_{n,m}$.
\end{proof}

The computation of $c^{(i)}_{\mu\nu}$, $d^{(i)}_{\mu\nu}$ is
straightforward and we state the result as follows.

For each permutation $\eta \in \Sigma_m$, there exists a unique
positive permutation braid (a positive braid that has a geometric
representative where every pair of strands crosses at most once, see
\cite{EM}) $\alpha_\eta \in \pi_1(\Y_{n,m}) \subset B_{n+m}$ in
which the last $m$ strands are permutated in the manner of $\eta$.
Set $\eta^\pm = \sgn\eta \cdot (\alpha_\eta)^{\pm 1}$.

For integers $1 \leq i \leq j \leq k \leq l \leq m$, let
$\theta_{i,j,k,l} \in \Sigma_m$ denote the permutation that sends
the sequence $i+1,i+2,\dots,l$ to
$$k+1,k+2,\dots,l, \;\;\; k,k-1,\dots,j+1, \;\;\; i+1,i+2,\dots,j.$$
Also set
$$\Theta_{j,k,l} = \Bigg\{ \eta\in\Sigma_m \Bigg|
 \begin{array}{l}
 \eta(i) = i, \; \forall i \neq j+1,j+2,\dots,l, \\
 \eta(j+1) < \eta(j+2) < \cdots < \eta(k), \\
 \eta(k+1) < \eta(k+2) < \cdots < \eta(l)
\end{array} \Bigg\}$$
and
$$\Theta_{j,k,l}^\pm = \sum_{\eta \in \Theta_{j,k,l}} \eta^\pm.$$

For integer $1 \leq i \leq n-1$ and elements
$\mu=(\mu_1,\dots,\mu_{n-1}), \nu=(\nu_1,\dots,\nu_{n-1}) \in
\E_{n,m}$, $c^{(i)}_{\mu\nu}$ and $d^{(i)}_{\mu\nu}$ do not vanish
if and only if $\mu_{i-1} \leq \nu_{i-1}$, $\mu_{i+1} \leq
\nu_{i+1}$ and $\mu_k = \nu_k$ for all $k \neq i-1,i,i+1$. In this
case,
\begin{eqnarray*}
 c^{(i)}_{\mu\nu} & = & (-1)^{\nu_i} \sigma_i
 \bigg( \prod_{k=u_{i+1}+1}^{v_i} A_{i,n+k} \bigg)
 \theta_{v_{i+1},v_{i+1},v_i,v_i}^+
 \Theta_{v_i,u_i,v_{i-1}}^+
 \Theta_{v_{i+2},u_{i+1},v_{i+1}}^+, \\
 d^{(i)}_{\mu\nu} & = & (-1)^{\nu_i}
 \theta_{u_{i+1},v_{i+1},v_i,u_i}^-
 \bigg( \prod_{k=v_{i+1}+1}^{u_i} A_{i,n+k} \bigg)^{-1}
 \Theta_{v_i,u_i,v_{i-1}}^+
 \Theta_{v_{i+2},u_{i+1},v_{i+1}}^+
 \sigma_i^{-1},
\end{eqnarray*}
where $u_j=\sum_{k=j}^{n-1}\mu_k$ and $v_j=\sum_{k=j}^{n-1}\nu_k$.

\subsection{The trace formula}

\begin{defn}
Let $\Gamma$ be a group, $\Z\Gamma$ its group ring, $\Gamma_c$ the
set of conjugacy classes, $\Z\Gamma_c$ the free abelian group
generated by $\Gamma_c$ and $\pi_\Gamma : \Z\Gamma \to \Z\Gamma_c$
the obvious projection. Suppose $\zeta$ is an endomorphism of a free
$\Z\Gamma$-module such that $\zeta(v_i) = \sum_{j=1}^k a_{ij} \cdot
v_j$ for a basis $\{ v_1,\dots,v_k \}$. The {\em trace} of $\zeta$
is defined as
$$\tr_\Gamma\zeta = \pi_\Gamma(\sum_{i=1}^k a_{ii}) \in \Z\Gamma_c.$$
It is straightforward to verify that the definition is independent
of the choice of the basis.
\end{defn}

Note that, under the basis $\E_{n,m}$, all matrix elements of
$\zeta_{n,m}(\beta)$ belong to $\Z\Gamma_{\beta,m}$. Therefore,
$\zeta_{n,m}(\beta)$ can be naturally regarded as an endomorphism of
the free $\Z\Gamma_{\beta,m}$-module generated by $\E_{n,m}$. In
this way, the notation  $\tr_{\Gamma_{\beta,m}} \zeta_{n,m}(\beta)$
in the following proposition makes sense.

Now we prove the main result of this section.

\begin{prop}\label{prop:trace}
For every nontrivial braid $\beta \in B_n$, we have
$$L_{\Gamma_{\beta,m}}(\SP^m\tf_\beta) =
  (-1)^m \tr_{\Gamma_{\beta,m}} \zeta_{n,m}(\beta) -
  \text{collapsible terms}.
$$
\end{prop}

\begin{proof}
Choose a word $\beta = \tau_1 \cdots \tau_k$ where
$\tau_1,\dots,\tau_k \in \{\sigma_1^{\pm1}, \dots,$
$\sigma_{n-1}^{\pm1}\}$. We put $\varphi_i = \phi_{j_i}$ if
$\tau_i=\sigma_{j_i}$ or $\varphi_i = \bar\phi_{j_i}$ if
$\tau_i=\sigma_{j_i}^{-1}$ . Then the embedding $g = \varphi_k
\cdots \varphi_1 : Y_n \to Y_n$ induces a map $\SP^mg:
\SP^mY_n\to\SP^mY_n$ stratified homotopic to $\SP^m\tf_\beta$. Hence
$L_{\Gamma_{\beta,m}}(\SP^m\tf_\beta) =
L_{\Gamma_{\beta,m}}(\SP^mg)$. It is immediate from the definitions
of $\phi_i$ and $\bar\phi_i$ that $\fix \DSP^mg \subset V$.

Note that the components of $\bigcup_{\mu\in\E_{n,m}} V_\mu \cap
\SP^mg^{-1}(V_\mu)$ are in one-one correspondence to the summands of
the last expression in the following equation.
\begin{eqnarray*}
  && (-1)^m \tr_{\Gamma_{\beta,m}} \zeta_{n,m}(\beta) \\
  &=& (-1)^m \tr_{\Gamma_{\beta,m}}
  \zeta_{n,m}(\tau_1) \cdots \zeta_{n,m}(\tau_k) \\
  &=& \sum_{\mu^0,\dots,\mu^k \in \E_{n,m}: \; \mu^0=\mu^k}
  \;\;\;\sum_{\eta^1,\dots,\eta^k\in\Sigma_m: \; W^{(\varphi_i)}_{\mu^{i-1}\mu^i\eta^i} \neq \emptyset} \\
  && (-1)^m (-1)^{\mu^1_{j_1} + \cdots + \mu^k_{j_k}}
  \sgn (\eta^1 \cdots \eta^k) \cdot
  [\tau_1\alpha^{(\varphi_1)}_{\mu^0\mu^1\eta^1} \cdots
  \tau_k\alpha^{(\varphi_k)}_{\mu^{k-1}\mu^k\eta^k}].
\end{eqnarray*}
Moreover, each of these components is homeomorphic to $\R^{2m}$ on
which $\SP^mg$ acts hyperbolically, hence gives rise to precisely
one fixed point of $\SP^mg$, either in $V$ or $\overline V \cap
\SP^{m,m-1}Y_n$. In the former case, the coordinate of the fixed
point corresponding to $\mu^0,\dots,\mu^k,\eta^1,\dots,\eta^k$ is
precisely
$$[\tau_1\alpha^{(\varphi_1)}_{\mu^0\mu^1\eta^1} \cdots
  \tau_k\alpha^{(\varphi_k)}_{\mu^{k-1}\mu^k\eta^k}].
$$
and, by Lemma \ref{lem:ind5}, the index is
$$(-1)^m (-1)^{\mu^1_{j_1} + \cdots + \mu^k_{j_k}} \sgn (\eta^1 \cdots \eta^k).$$
In the latter case, the corresponding summand is always collapsible.
Therefore, from Lemma \ref{lem:collapse} the proposition follows.
\end{proof}

\begin{rem}
In fact, the configuration space $\Y_{n,m}$ has the homotopy type of
a compact $m$-complex and the trace $(-1)^m \tr_{\Gamma_{\beta,m}}
\zeta_{n,m}(\beta)$ is nothing but the generalized Lefschetz number
of a self map of the complex induced by $\DSP^mf_\beta$. In this
sense, the collapsible and peripheral terms in the trace both arise
from the compactification issue.
\end{rem}

\section{Proof of main theorem}\label{sec:conjugacy}

According to Proposition \ref{prop:Lefschetz} and Proposition
\ref{prop:trace}, the $(n+m)$-strand forced extensions of a
nontrivial braid $\beta\in B_n$ are precisely those non-collapsible,
non-peripheral terms in the trace $\tr_{\Gamma_{\beta,m}}
\zeta_{n,m}(\beta)$. The following proposition states that these
terms do not cancel in $\tr_{B_{n+m}} \zeta_{n,m}(\beta)$, hence
eventually establishes Theorem \ref{thm:main}.

\begin{prop}\label{prop:conj}
Let $\beta \in B_n$ be a nontrivial braid and suppose two extensions
$\beta', \beta''\in \beta\cdot\pi_1(\Y_{n,m})$ of $\beta$ are
conjugate in $B_{n+m}$. If $\beta'$ is forced by $\beta$, then
$[\beta']$ and $[\beta'']$ have the same coefficient in
$L_{\Gamma_{\beta,m}}(\SP^m\tf_\beta)$.
\end{prop}

\begin{proof}
Assume $f_\beta=\tf_\beta|_{\intr Y_n}: \R^2 \setminus P\to\R^2
\setminus P$ is a minimal representative (in the sense of Theorem
\ref{thm:min2}) in its isotopy class, and assume the term $[\beta']$
in $L_{\Gamma_{\beta,m}}(\SP^m\tf_\beta)$ is the coordinate of the
fixed point of $\DSP^m\tf_\beta$ corresponding to an
$\tf_\beta$-invariant set $Q\subset\intr Y_n=\R^2 \setminus P$. We
extend $f_\beta$ to a homeomorphism $\phi:\R^2\to\R^2$. Suppose the
puncture point set $P$ splits into a disjoint union of periodic
orbits $c_1 \cup \cdots \cup c_s$ of $\phi$ and suppose $Q$ splits
into a disjoint union of periodic orbits $d_1 \cup \cdots \cup d_t$
of $\phi$.

The conjugation between $\beta'$ and $\beta^{\prime\prime}$ in
$B_{n+m}$ gives rise to a homeomorphism $\psi : \R^2 \to \R^2$,
which preserves the set $P \cup Q = c_1 \cup \dots \cup c_s \cup d_1
\cup \dots \cup d_t$. Put $\phi' = \psi\phi\psi^{-1}$. Since
$\beta''$ restricts to $\beta$ on the first $n$ strands,
$\phi'|_{\R^2\setminus P}$ is isotopic to $f_\beta$. Further,
$\phi'|_{\R^2\setminus P}$ is also a minimal representative in its
isotopy class.

By conjugating $\beta^{\prime\prime}$ in $\Gamma_{\beta,m}$ if
necessary, we may assume $d_1,\dots,d_t$ are periodic orbits of
$\phi'$. Let $m_i$ be the period of $d_i$. Then
$$\ind(\phi^{m_i},\psi^j(d_i)) = \ind(\phi^{\prime m_i},\psi^{j+1}(d_i))$$
for all $j$ and
$$\ind(\phi^{m_i},\psi^j(d_i)) = \ind(\phi^{\prime m_i},\psi^j(d_i))$$
provided that $\psi^j(d_i) \in \{ c_1,\dots,c_s \}$.

Let $n_i$ be the maximum positive number such that
$$\psi(d_i),\dots,\psi^{n_i-1}(d_i) \in \{ c_1,\dots,c_s \}.$$
Then
$$\{ \psi^{n_1}(d_1),\dots,\psi^{n_t}(d_t) \} = \{ d_1,\dots,d_t \},$$
and, by induction,
$$\ind(\phi^{m_i},d_i) = \ind(\phi^{\prime m_i},\psi^{n_i}(d_i)).$$
On the other hand, by Lemma \ref{lem:ind2} and Lemma \ref{lem:ind3},
we have
\begin{eqnarray*}
  \ind(\DSP^m\tf_\beta,[\beta'])
  &=& \prod_{i=1}^t \ind(\phi^{m_i},d_i) / m_i, \\
  \ind(\DSP^m\tf_\beta,[\beta^{\prime\prime}])
  &=& \prod_{i=1}^t \ind(\phi^{\prime m_i},d_i) / m_i.
\end{eqnarray*}
Therefore, these two indices are identical.
\end{proof}

\section{Algorithms and examples}\label{sec:algorithm}

Thanks to Theorem \ref{thm:main}, the computation of the
$(n+m)$-strand forced extensions of a given braid $\beta\in B_n$ may
proceed as follows.
\begin{enumerate}
\setlength{\itemsep}{-.5ex}
\item
By means of the representation $\zeta_{n,m}$, compute an initial
formal sum for the trace $\tr_{B_{n+m}} \zeta_{n,m}(\beta)$.

\item
Merge conjugate terms in the formal sum by solving conjugacy problem
in $B_{n+m}$.

\item
Identify collapsible terms and peripheral terms by computing
reducing curves and drop them off.

\item
Return the nonzero terms remained after cancelation.
\end{enumerate}

In the procedure described above, one has to deal with two
algorithmic problems: conjugacy problem in the braid group $B_{n+m}$
and computation of reducing curves. Fortunately, there have been
effective algorithms for both tasks.

For the conjugacy problem, we refer the reader to a very efficient
algorithm due to Gebhardt \cite{Gebhardt}. See also \cite{BGG,
Zheng2} for improvements on this direction.

As to the second problem, one solution is a braid algorithm due to
Bernardete-Nitecki-Gutierrez \cite{BNG}. It can be improved
significantly if one computes the ultra summit set \cite{Gebhardt}
or its variant \cite{Zheng2} instead of the super summit set (see
the references for details). An alternative solution is given by
Bestvina-Handel \cite{BH}, which is also applicable for general
surface homeomorphisms but apparently less efficient, because it
involves a computation of train-track maps.

\smallskip

At the present time, we are not able to talk much about the
computational complexity of the above procedure, partly because the
topic of braid algorithms is a fairly new one and many questions
still remain open. Nevertheless, the bulk part of running time is
evidently spent in the second step. Hence it is a major issue to
control the number of terms written down in the first step.

A braid is called {\em cyclic} if it induces a cyclic permutation on
the end points of its strands. We call an extension of a braid {\em
elementary} if it is obtained by appending a single cyclic braid.

Elementary forced extensions are the main concern of the braid
forcing problem. Observe that the elementary extensions only
constitute a small fraction of the terms in
$L_{\Gamma_{\beta,m}}(\SP^m\tf_\beta)$. Hence in the first step of
the above procedure we may drop off all non-elementary extensions of
$\beta$ to save considerably on running time.

\smallskip

As another example of shortcut to facilitate the computation, when
$\beta$ is a pseudo-Anosov braid (the most significant case in
dynamics), identification of collapsible and peripheral terms may be
reduced by the following proposition to the reducibility problem of
braids, for which a polynomial solution (for a fixed number of
strands) has been claimed recently by Ko-Lee \cite{KL}.

\begin{prop}\label{prop:pA}
An extension $\beta'$ of a pseudo-Anosov $\beta$ is collapsible or
peripheral relative to $\beta$ if and only if $\beta'$ is a
reducible braid.
\end{prop}
\begin{proof}
Let $\phi$ be a Thurston representative determined by $\beta'$. If
$\beta'$ is not reducible, there are no reducing curves of $\phi$,
hence $\beta'$ is neither collapsible nor peripheral relative to
$\beta$. Conversely, If $\beta'$ is reducible, each reducing curve
of $\phi$ must enclose either at most one of or all of the punctures
corresponding to $\beta$, because $\beta$ is pseudo-Anosov. Hence
$\beta'$ is either collapsible or peripheral relative to $\beta$.
\end{proof}

Below we conclude this paper by presenting some examples.

\begin{exam}
Under the basis $\E_{3,2} = \{ (2,0),(1,1),(0,2) \}$ the
representation $\zeta_{3,2}$ of $B_3$ is given by the matrices (cf.
the matrices of $\xi_{3,2}$ from \cite{Zheng1})

$$\begin{array}{ccl}
  \zeta_{3,2}(\sigma_1) &=& \sigma_1 \cdot \left( \begin{array}{ccc}
  -A_{14}A_{15}\sigma_4 & -A_{14}A_{15} & A_{14}A_{15} \\
  0                     & -A_{15}       & A_{15}(1-\sigma_4) \\
  0                     & 0             & 1
  \end{array} \right),\\
\\
  \zeta_{3,2}(\sigma_1^{-1}) &=& \left( \begin{array} {ccc}
  -\sigma_4^{-1}A_{15}^{-1}A_{14}^{-1} & \sigma_4^{-1}A_{15}^{-1} & 1 \\
  0                                    & -A_{15}^{-1}             & 1-\sigma_4 \\
  0                                    & 0                        & 1
  \end{array} \right) \cdot \sigma_1^{-1},\\
\\
  \zeta_{3,2}(\sigma_2) &=& \sigma_2 \cdot \left( \begin{array}{ccc}
  1          & 0       & 0 \\
  1-\sigma_4 & -A_{24} & 0 \\
  1          & -A_{24} & -A_{24}A_{25}\sigma_4
  \end{array} \right),\\
\\
  \zeta_{3,2}(\sigma_2^{-1}) &=& \left( \begin{array} {ccc}
  1                       & 0                                   & 0 \\
  A_{24}^{-1}(1-\sigma_4) & -A_{24}^{-1}                        & 0 \\
  A_{25}^{-1}A_{24}^{-1}  & \sigma_4^{-1}A_{25}^{-1}A_{24}^{-1} & -\sigma_4^{-1}A_{25}^{-1}A_{24}^{-1}
  \end{array} \right) \cdot \sigma_2^{-1}.
\end{array}
$$

For reader's convenience, we illustrate by figures how to obtain the
equality
$$(1,1)\cdot\zeta_{3,2}(\sigma_1) = -\sigma_1A_{15}\cdot(1,1) + \sigma_1A_{15}(1-\sigma_4)\cdot(0,2).$$
See Figure \ref{fig:fig31}. The set $V_{(1,1)}$ consists of those
points $[x_1,x_2]$ with $x_1,x_2\in Y_3$ positioned as in the top
left figure. Note that the set $V_{(1,1)} \cap
\phi_1^{-1}(V_{(2,0)})$ is empty; $V_{(1,1)} \cap
\phi_1^{-1}(V_{(1,1)})$ has one component, illustrated by the top
right figure; and $V_{(1,1)} \cap \phi_1^{-1}(V_{(0,2)})$ has two
components, illustrated by the bottom two figures. From the last
three figures one reads out $A_{15}$, $A_{15}$ and $A_{15}\sigma_4$,
respectively. Together with the contribution of signs, they are
assembled to give the above equality.

\begin{figure}[h]
\centering
\includegraphics[scale=.6]{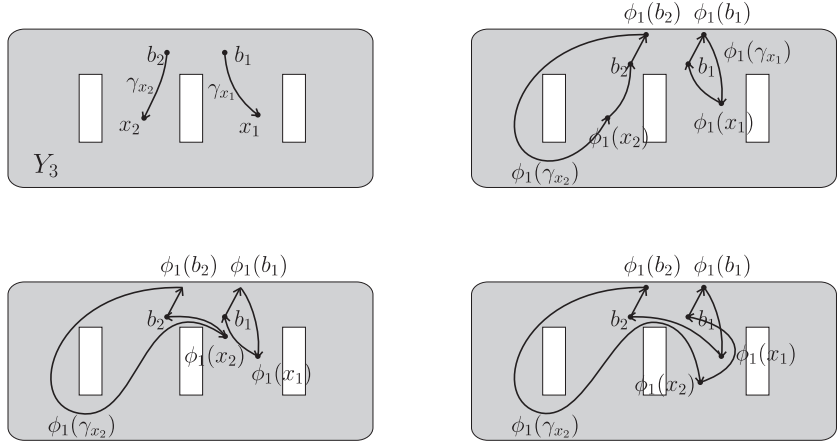}
\caption{ figures for computing $(1,1)\cdot\zeta_{3,2}(\sigma_1)$ }
\label{fig:fig31}
\end{figure}
\end{exam}

\begin{exam}
Under the basis $\E_{n,1}$ the representation $\zeta_{n,1}$ of $B_n$
is given by the matrices
$$\zeta_{n,1}(\sigma_i) = \sigma_i \cdot
  \left( \begin{array}{ccccc}
    I_{i-2} \\
    & 1         & 0          & 0 \\
    & A_{i,n+1} & -A_{i,n+1} & 1 \\
    & 0         & 0          & 1 \\
    & & & & I_{n-i-2}
  \end{array} \right).
$$
Note that if we replace $\sigma_i$ by $1$ and replace $A_{i,n+1}$ by
a number $a$, the representation specializes to the reduced Burau
representation (ref. \cite{Birman})
$$\sigma_i \mapsto
  \left( \begin{array}{ccccc}
    I_{i-2} \\
    & 1         & 0          & 0 \\
    & a         & -a         & 1 \\
    & 0         & 0          & 1 \\
    & & & & I_{n-i-2}
  \end{array} \right).
$$
\end{exam}

\begin{exam}
For the simplest pseudo-Anosov braid $\beta =
\sigma_1\sigma_2^{-1}$, we have
\begin{eqnarray*}
  \tr_{B_5} \zeta_{3,2}(\beta)
  &=& \tr_{B_5} \zeta_{3,2}(\sigma_1) \zeta_{3,2}(\sigma_2^{-1}) \\
  &=& [\sigma_1\cdot(
  -A_{14}A_{15}\sigma_4
  -A_{14}A_{15}A_{24}^{-1}(1-\sigma_4) \\ &&
  +A_{14}A_{15}A_{25}^{-1}A_{24}^{-1}
  +A_{15}A_{24}^{-1} \\ &&
  +A_{15}(1-\sigma_4)\sigma_4^{-1}A_{25}^{-1}A_{24}^{-1}
  -\sigma_4^{-1}A_{25}^{-1}A_{24}^{-1}
  )\cdot\sigma_2^{-1} ] \\
  &=& [
  \beta
  -\beta A_{35}^{-1}A_{34}^{-1}\sigma_4^{-1}
  -\beta A_{14}A_{15}\sigma_4
  -\beta A_{34}^{-1}
  -\beta A_{15} \\ &&
  +\beta A_{34}^{-1}\sigma_4^{-1}
  +\beta A_{15}\sigma_4
  +\beta A_{15}A_{34}^{-1}
  ].
\end{eqnarray*}

\begin{figure}[h]
    \psfrag{1}[]{$\beta$}
    \psfrag{2}[]{$\beta A_{35}^{-1}A_{34}^{-1}\sigma_4^{-1}$}
    \psfrag{3}[]{$\beta A_{14}A_{15}\sigma_4$}
    \psfrag{4}[]{$\beta A_{34}^{-1}$}
    \psfrag{5}[]{$\beta A_{15}$}
    \psfrag{6}[]{$\beta A_{34}^{-1}\sigma_4^{-1}$}
    \psfrag{7}[]{$\beta A_{15}\sigma_4$}
    \psfrag{8}[]{$\beta A_{15}A_{34}^{-1}$}
\centering
\includegraphics[scale=.6]{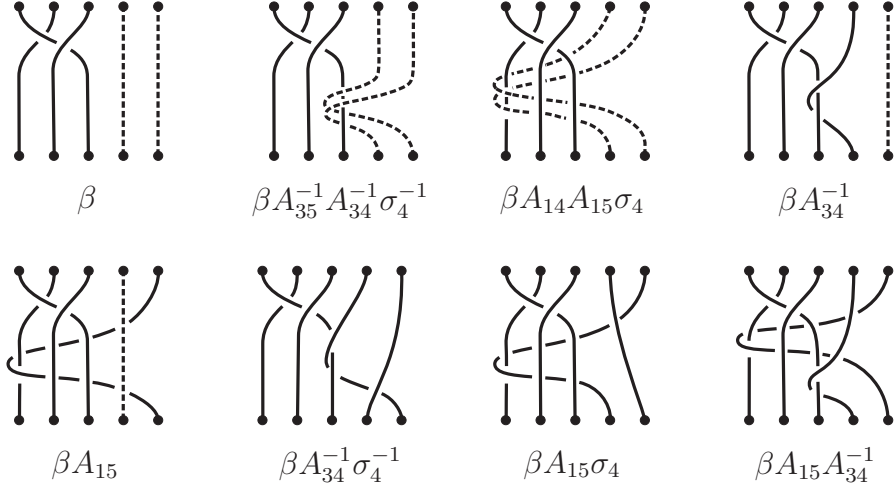}
\caption{ the braids appearing in $\tr_{B_5} \zeta_{3,2}(\beta)$ for
$\beta=\sigma_1\sigma_2^{-1}$ } \label{fig:fig32}
\end{figure}

See Figure \ref{fig:fig32}, in which the collapsible or peripheral
strands are depicted as dotted lines. Clearly, the first five braids
in the figure are reducible. An algorithmic test shows the last
three are pseudo-Anosov. It follows from Proposition \ref{prop:pA}
that precisely the last three terms in $\tr_{B_5}
\zeta_{3,2}(\beta)$ are neither collapsible nor peripheral.
Therefore, up to conjugacy, there are a total of three $5$-strand
forced extensions of $\beta$: $\beta A_{34}^{-1}\sigma_4^{-1}$,
$\beta A_{15}\sigma_4$ and $\beta A_{15}A_{34}^{-1}$.
\end{exam}

\begin{exam}\label{exam:4}
Suppose $\beta = \sigma_1 \cdots \sigma_{n_1} \sigma_{n_1+1}^{-1}
\cdots \sigma_{n_1+n_2}^{-1} \in B_n$ where $n_1, n_2 \geq 2$ and
$n=n_1+n_2+1$. For $2 \leq m \leq \min(n_1,n_2)$, we have
\begin{align*}
  & \tr_{B_{n+1}} \zeta_{n,1}(\beta^m) =
  [\beta^m - (\beta A_{1,n+1})^m - (\beta A_{n,n+1}^{-1})^m], \\
  & \tr_{B_{n+m}} \zeta_{n,m}(\beta) =
  [\beta(1-A_{1,n+2})(1-A_{n,n+1}^{-1})].
\end{align*}
Either of the above formulae implies that the (pseudo-Anosov) cyclic
braid $\beta$ forces no $m$-strand cyclic braid (see \cite[Theorem
7]{Guaschi2} for the case $n_2=m=2$). This contrasts sharply to
Guaschi's theorem \cite{Guaschi1} which asserts that a pseudo-Anosov
braid on three or four strands forces at least one $m$-strand cyclic
braid for every $m\ge1$.

Without loss of generality, we sketch the computation of the above
formulae for $n_1=n_2=m=2$. First, we translate the braid $\beta =
\sigma_1 \sigma_2 \sigma_3^{-1} \sigma_4^{-1} \in B_5$ into the self
embedding $\phi: Y_5\to Y_5$ depicted in Figure \ref{fig:fig33}.

\begin{figure}[h]
\centering
\includegraphics[scale=.6]{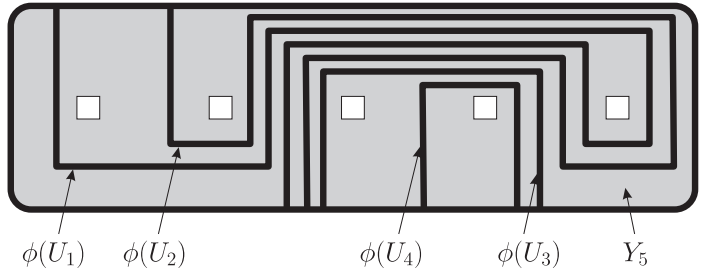}
\caption{ the map $\phi: Y_5\to Y_5$ representing $\beta = \sigma_1
\sigma_2 \sigma_3^{-1} \sigma_4^{-1}$ } \label{fig:fig33}
\end{figure}

Keep the notations of Section \ref{sec:formula}. There are a total
of 11 components in $V\cap\phi^{-1}(V)$, so there should be the same
number of nonzero terms in the matrix $\zeta_{5,1}(\beta)$ under the
standard basis $\E_{5,1}$
$$\zeta_{5,1}(\beta) = \beta\cdot\begin{pmatrix}
  0 & A_{1,6}A_{2,6}(-1+A_{5,6}^{-1}) & 0 & A_{1,6}A_{2,6}A_{5,6}^{-1} \\
  1 & A_{2,6}(-1+A_{5,6}^{-1}) & 0 & A_{2,6}A_{5,6}^{-1} \\
  0 & A_{5,6}^{-1} & 0 & -A_{5,6}^{-1} \\
  0 & 0 & A_{5,6}^{-1} & -A_{5,6}^{-1} \\
  \end{pmatrix}.
$$
Since the matrix is almost upper triangular (this is quite evident
for larger $n_1,n_2$), the following trace can be computed without
much difficulty.
\begin{eqnarray*}
  \tr_{B_6} \zeta_{5,1}(\beta^2)
  &=& [\beta A_{1,6}A_{2,6}(-1+A_{5,6}^{-1})\beta
    +\beta^2 A_{1,6}A_{2,6}(-1+A_{5,6}^{-1})
    \\&&
    +(\beta A_{2,6}(-1+A_{5,6}^{-1}))^2
    -(\beta A_{5,6}^{-1})^2
    -(\beta A_{5,6}^{-1})^2
    +(\beta A_{5,6}^{-1})^2
  ] \\
  &=& [\beta^2 - (\beta A_{1,6})^2 - (\beta A_{5,6}^{-1})^2].
\end{eqnarray*}
In the equality, we used the identities $\beta A_{1,6}=A_{2,6}\beta$
and $[\beta A_{2,6}A_{5,6}^{-1}]=[\beta]$.

Next, we compute the second formula. The matrix $\zeta_{5,2}(\beta)$
is a $10\times10$ one, but we are only concerned with its diagonal
part. Notice that there are totally 10 components in
$\bigcup_{\mu\in\E_{5,2}} V_\mu \cap \phi^{-1}(V_\mu)$. The trace is
computed as
\begin{eqnarray*}
  \tr_{B_7} \zeta_{5,2}(\beta)
  &=& [-\beta \sigma_6A_{1,6}A_{2,6}(-1+A_{5,6}^{-1})
    -\beta A_{2,6}A_{2,7}\sigma_6
    \\&&
    -\beta A_{2,6}A_{2,7}A_{5,6}^{-1}(1-\sigma_6)
    +\beta A_{2,6}A_{2,7}(A_{5,6}A_{5,7})^{-1}
    \\&&
    -\beta A_{2,7}(-1+A_{5,7}^{-1})A_{5,6}^{-1}
    +\beta(A_{5,6}A_{5,7})^{-1}\sigma_6^{-1}
    -\beta(A_{5,6}A_{5,7})^{-1}\sigma_6^{-1}
  ].
\end{eqnarray*}
On the right hand side, the six terms containing $\sigma_6^{\pm1}$
cancel pairwise, and we get
\begin{eqnarray*}
  &=& [-\beta A_{2,6}A_{2,7}A_{5,6}^{-1}
    +\beta A_{2,6}A_{2,7}(A_{5,6}A_{5,7})^{-1}
    -\beta A_{2,7}(-1+A_{5,7}^{-1})A_{5,6}^{-1}
  ] \\
  &=& [\beta(1-A_{1,7})(1-A_{5,6}^{-1})].
\end{eqnarray*}
\end{exam}


\begin{thebibliography}{99}
\setlength{\itemsep}{-.5ex}

\bibitem{BNG} D. Bernardete, Z. Nitecki and M. Gutierrez,
    Braids and the Nielsen-Thurston classification,
    J. Knot Theory Ramif. 4 (1995) 549--618.
\bibitem{BH} M. Bestvina and M. Handel,
    Train tracks for surface homeomorphisms,
    Topology 34 (1995) 109--140.
\bibitem{Birman} J. S. Birman,
    Braids, Links, and Mapping Class Groups,
    Ann. Math. Stud. 82, Princeton Univ. Press, Princeton, 1974.
\bibitem{BGG} J. S. Birman, V. Gebhardt and J. Gonz\'alez-Meneses,
    Conjugacy in Garside groups I, II, III,
    arXiv:math.GT/0605230, 0606652, 0609616.
\bibitem{Boyland1} P. Boyland,
    Braid types and a topological method of proving positive entropy,
    Boston University, 1984, preprint.
\bibitem{Boyland2} P. Boyland,
    Isotopy stability of dynamics on surfaces,
    in: Geometry and Topology in Dynamics,
    Contemp. Math. 246, Amer. Math. Soc., Providence, 1999, pp. 17--46.
\bibitem{CH1} A. de Carvalho and T. Hall,
    The forcing relation for horseshoe braid types,
    Exp. Math. 11 (2002) 271--288.
\bibitem{CH2} A. de Carvalho and T. Hall,
    Braid forcing and star-shaped train tracks,
    Topology 43 (2004) 247--287.
\bibitem{EM} E. A. Elrifai and H. R. Morton,
    Algorithms for positive braids,
    Quart. J. Math. Oxford 45 (1994) 479--497.
\bibitem{FLP} A. Fathi, F. Lauderbach and V. Po\'enaru,
    Travaux de Thurston sur les surfaces,
    Ast\'erique 66-67 (1979).
\bibitem{Gebhardt} V. Gebhardt,
    A new approach to the conjugacy problem in Garside groups,
    J. Algebra 292(1) (2005) 282--302.
\bibitem{Guaschi1} J. Guaschi,
    Pseudo-Anosov braid types of the disc or sphere of low cardinality imply all periods,
    J. London Math. Soc. (2) 50 (1994) 594--608.
\bibitem{Guaschi2} J. Guaschi,
    Nielsen theory, braids and fixed points of surface homeomorphisms,
    Topol. Appl. 117 (2002) 199--230.
\bibitem{Handel} M. Handel,
    The forcing partial order on the three times punctured disk,
    Ergod. Theory Dynam. Sys. 17 (1997) 593--610.
\bibitem{Jiang1} B. Jiang,
    Lectures on Nielsen Fixed Point Theory,
    Contemp. Math. 14, Amer. Math. Soc., Providence, 1983.
\bibitem{Jiang2} B. Jiang,
    Estimation of the number of periodic orbits,
    Pacific J. Math. 172 (1996) 151--185.
\bibitem{JG} B. Jiang and J. Guo,
    Fixed points of surface diffeomorphisms,
    Pacific J. Math. 160 (1993) 67--89.
\bibitem{JZZ} B. Jiang, X. Zhao and H. Zheng,
    On fixed points of stratified maps,
    preprint.
\bibitem{KL} K. H. Ko and J. W. Lee,
    A polynomial-time solution to the reducibility problem,
    arXiv:math.GT/0610746.
\bibitem{Kolev} B, Kolev,
    Periodic points of period 3 in the disc,
    Nonlinearity 7 (1994) 1067--1072.
\bibitem{LY} T. Y. Li and J. A. Yorke,
    Period three implies chaos,
    Amer. Math. Monthly 82 (1975) 985--992.
\bibitem{Matsuoka} T. Matsuoka,
    The number and linking of periodic solutions of periodic systems,
    Invent. Math. 70 (1983) 319--340.
\bibitem{Sharkovskii} A. N. Sharkovskii,
    Coexistence of cycles of a continuous map of a line to itself,
    Ukr. Math. Z. 16 (1964) 61--71.
\bibitem{Thurston} W. P. Thurston,
    On the geometry and dynamics of diffeomorphisms of surfaces,
    Bull. Amer. Math. Soc. 19 (1988) 417--431.
\bibitem{Zheng1} H. Zheng,
    A reflexive representation of braid groups,
    J. Knot Theory Ramif. 14 (2005) 467--477.
\bibitem{Zheng2} H. Zheng,
    General cycling operations in Garside groups,
    arXiv:math.GT/ 0605741.

\end{thebibliography}
\end{document}